\documentclass[12pt]{article}
\usepackage{amsmath}
\usepackage{amsfonts}
\usepackage{amssymb}
\usepackage{graphicx}%
\providecommand{\U}[1]{\protect \rule{.1in}{.1in}}
\newtheorem{theorem}{Theorem}

\newtheorem{corollary}[theorem]{Corollary}

\newtheorem{definition}[theorem]{Definition}

\newtheorem{lemma}[theorem]{Lemma}

\newtheorem{remark}[theorem]{Remark}

\newenvironment{proof}[1][Proof]{\noindent \textbf{#1.} }{\  $\Box$}

\title{Choquet expectations and $g$-expectations with multi-dimensional
Brownian motion}

\author{Mingshang Hu\footnote{Email: humingshang\symbol{64}gmail.com}\\
        \footnotesize{
        School of Mathematics, Shandong University, Jinan, \small{250100}, China
 } }

\begin{document}
\maketitle

\noindent \textbf{Abstract. }We prove that a $g$-expectation is a
Choquet expectation if and only if $g$ is independent of $y$ and is
linear in $z$, i.e., classical linear expectation, without the
assumptions that the deterministic generator $g$ is continuous in
$t$ and the dimension of the Brownian motion is one.

\  \  \

\noindent \textbf{Keywords: }BSDE, $g$-expectation, conditional
$g$-expectation, capacity, Choquet expectation, comonotonic
additivity.

\section{Introduction}

Choquet \cite{Choq} introduced the notion of Choquet expectations
via capacities in 1953. Peng \cite{P2} introduced the notions of
$g$-expectations and conditional $g$-expectations via a class of
backward stochastic differential equations (BSDEs for short) in
1997. These two types of nonlinear mathematical expectations have
their own characteristics. For example, Choquet expectations are
comonotonic additivity, $g$-expectations and conditional
$g$-expectations are consistent. In Chen et al. \cite{CCD}, the
authors studied an interesting problem:

If a $g$-expectation is a Choquet expectation, can we find the form
of the generator $g$?

Under the assumptions that the deterministic generator $g$ is
continuous in $t$ and the dimension of the Brownian motion is one,
Chen et al. \cite{CCD} proved that a $g$-expectation is a Choquet
expectation if and only if $g$ is independent of $y$ and is linear
in $z$. For the case that the dimension of the Brownian motion is
greater than one, the main difficulty is to find the form of the
generator $g$. Unfortunately, this problem is not a simple extension
of the one dimensional case. Take a $2$-dimensional Brownian motion
$W_{t}=(W_t^1,W_t^2)$ for example, $W_t^1$ and $W_t^2$ are not
comonotonic. This prevents us from using the method in Chen et al.
\cite{CCD} directly. To overcome this defect, we consider
comonotonic indicator functions and use a property of BSDE.
Furthermore, our method does not need the continuous assumption on
$g$.

This paper is organized as follows: In Section 2, we recall some facts about
$g$-expectations and Choquet expectations. In Section 3, we state and prove
our main result.

\section{Preliminaries}

Fix $T>0$, let $(W_{t})_{0\leq t\leq T}$ be a $d$-dimensional standard
Brownian motion defined on a completed probability space $(\Omega
,\mathcal{F},P)$ and $(\mathcal{F}_{t})_{0\leq t\leq T}$ be the natural
filtration generated by this Brownian motion. For $x=(x_{1},\ldots,x_{d})$,
$y=(y_{1},\ldots,y_{d})\in \mathbb{R}^{d}$, $|x|:=\sqrt{\sum_{i=1}^{d}%
|x_{i}|^{2}}$, $x\cdot y:=\sum_{i=1}^{d}x_{i}y_{i}$. We denote by
$L^{2}(\mathcal{F}_{t})$ the set of all square integrable $\mathcal{F}_{t}%
$-measurable random variables and $L^{2}(0,T;\mathbb{R}^{n})$ the space of all
$\mathcal{F}_{t}$-adapted, $\mathbb{R}^{n}$-valued processes $(v_{t}%
)_{t\in \lbrack0,T]}$ with $E\int_{0}^{T}|v_{t}|^{2}dt<\infty$.

Let us consider a deterministic function $g:[0,T]\times \mathbb{R\times R}%
^{d}\rightarrow \mathbb{R}$, which will be in the following the generator of
the BSDE. For the function $g$, we will use the following assumptions:

\begin{description}
\item[(H1)] For each $(y,z)\in \mathbb{R\times R}^{d}$, $t\rightarrow g(t,y,z)$
is measurable.

\item[(H1$^{\prime}$)] For each $(y,z)\in \mathbb{R\times R}^{d}$,
$t\rightarrow g(t,y,z)$ is continuous.

\item[(H2)] There exists a constant $K\geq0$ such that
\[
|g(t,y,z)-g(t,y^{\prime},z^{\prime})|\leq K(|y-y^{\prime}|+|z-z^{\prime
}|),\  t\in \lbrack0,T],y,y^{\prime}\in \mathbb{R},z,z^{\prime}%
\in \mathbb{ R}^{d}.
\]

\item[(H3)] $g(t,y,0)\equiv0$ for each $(t,y)\in \lbrack0,T]\times \mathbb{R}.$
\end{description}

Let $g$ satisfy (H1)-(H3). Then for each $\xi \in L^{2}(\mathcal{F}_{T})$, the
BSDE%
\begin{equation}\label{e1}
y_{t}=\xi+\int_{t}^{T}g(s,y_{s},z_{s})ds-\int_{t}^{T}z_{s}\cdot dW_{s},\ 0\leq
t\leq T,
\end{equation}
has a unique solution $(y_{t,}z_{t})_{t\in \lbrack0,T]}\in L^{2}(0,T;\mathbb{R}%
)\times L^{2}(0,T;\mathbb{R}^{d})$ (see Pardoux and Peng \cite{PP}),
which depends on the generator $g$ and terminal value $\xi.$

The following standard estimate for BSDEs can be found in \cite{EPQ,
P1, BCHMP}.

\begin{lemma}\label{le1}
Suppose $g$ satisfies (H1)-(H3). For each $\xi^{1},\xi^{2}\in L^{2}%
(\mathcal{F}_{T})$, let $(y_{t,}^{i}z_{t}^{i})_{t\in \lbrack0,T]}$
be the solution of BSDE (\ref{e1}) corresponding to the generator
$g$ and terminal value
$\xi^{i}$ with $i=1,2.$ Then there exists a constant $C>0$ such that%
\[
E[\sup_{t\leq s\leq T}|y_{s}^{1}-y_{s}^{2}|^{2}|\mathcal{F}_{t}]+E[\int
_{t}^{T}|z_{s}^{1}-z_{s}^{2}|^{2}ds|\mathcal{F}_{t}]\leq CE[|\xi^{1}-\xi
^{2}|^{2}|\mathcal{F}_{t}].
\]

\end{lemma}

Using the solution of BSDE (\ref{e1}), Peng \cite{P2} proposed the
following notions:

\begin{definition}\label{de2}
Suppose $g$ satisfies (H1)-(H3). For each $\xi \in L^{2}(\mathcal{F}_{T}),$ let
$(y_{t,}z_{t})_{t\in \lbrack0,T]}$ be the solution of BSDE (\ref{e1}), define%
\[
\mathcal{E}_{g}[\xi]:=y_{0};\qquad \mathcal{E}_{g}[\xi|\mathcal{F}_{t}%
]:=y_{t}\quad \mbox{for each}\ t\in \lbrack0,T].
\]
$\mathcal{E}_{g}[\xi]$ is called the $g$-expectation of $\xi$ and
$\mathcal{E}_{g}[\xi|\mathcal{F}_{t}]$ is called the conditional
$g$-expectation of $\xi$ with respect to $\mathcal{F}_{t}$.
\end{definition}

We now recall the notions of capacity and Choquet expectation. A
capacity is a set function $V:\mathcal{F}_{T}\mapsto \lbrack0,1]$
satisfying: (i) $V(\emptyset)=0,$ $V(\Omega)=1;$ (ii) $V(A)\leq
V(B)$ for each $A\subset B.$ The corresponding Choquet expectation
(see \cite{Choq}) is defined as
follows:%
\[
\mathcal{C}[\xi]:=\int_{-\infty}^{0}[V(\xi \geq t)-1]dt+\int_{0}^{\infty}%
V(\xi \geq t)dt\quad \mbox{for each}\  \xi \in L^{2}(\mathcal{F}_{T}).
\]
Two random variables $\xi$ and $\eta$ are called comonotonic if%
\[
\lbrack \xi(\omega)-\xi(\omega^{\prime})][\eta(\omega)-\eta(\omega^{\prime
})]\geq0\quad \mbox{for each}\  \omega,\omega^{\prime}\in \Omega.
\]
Now, we list some properties of Choquet expectations (see
\cite{Choq, S1, De, Den}).

\begin{description}
\item[(1)] Monotonicity: If $\xi \geq \eta,$ then $\mathcal{C}[\xi
]\geq \mathcal{C}[\eta].$

\item[(2)] Positive homogeneity: If $\lambda \geq0,$ then $\mathcal{C}%
[\lambda \xi]=\lambda \mathcal{C}[\xi].$

\item[(3)] Translation invariance: If $c\in \mathbb{R},$ then $\mathcal{C}%
[\xi+c]=\mathcal{C}[\xi]+c.$

\item[(4)] Comonotonic additivity: If $\xi$ and $\eta$ are comonotonic, then
$\mathcal{C}[\xi+\eta]=\mathcal{C}[\xi]+\mathcal{C}[\eta].$
\end{description}

Let $g$ satisfy (H1)-(H3), define%
\[
P_{g}(A):=\mathcal{E}_{g}[I_{A}]\quad \mbox{for each}\ A\in \mathcal{F}_{T}.
\]
$P_{g}(A)$ is called the $g$-probability of $A.$ Obviously, $P_{g}(\cdot)$ is
a capacity. The corresponding Choquet expectation is denoted by $\mathcal{C}%
_{g}$. It is easy to check that $\mathcal{C}_{g}[I_{A}]=\mathcal{E}_{g}%
[I_{A}]$ for each $A\in \mathcal{F}_{T}$. Furthermore, $\mathcal{C}_{g}%
[\xi]<\infty$ for each $\xi \in L^{2}(\mathcal{F}_{T})$ (see
\cite{HHC}).

The following result can be found in \cite{CCD}.

\begin{lemma}\label{le3}
Suppose that $d=1$ and $g$ satisfies (H1$^{\prime}$), (H2) and (H3). Then
$\mathcal{E}_{g}[\xi]=\mathcal{C}_{g}[\xi]$ for each $\xi \in L^{2}%
(\mathcal{F}_{T})$ if and only if $g$ is independent of $y$ and is linear in
$z$, i.e., $g(t,z)=g(t,1)z$.
\end{lemma}

\section{Main result}

Let $\{e_{1},e_{2},\ldots,e_{d}\}$ denote the standard basis of $\mathbb{R}%
^{d}$. Now we give the main result.

\begin{theorem}\label{th4}
Suppose $g$ satisfies (H1)-(H3). Then $\mathcal{E}_{g}[\xi]=\mathcal{C}%
_{g}[\xi]$ for each $\xi \in L^{2}(\mathcal{F}_{T})$ if and only if $g$ is
independent of $y$ and is linear in $z$, i.e., $g(t,z)=\sum_{i=1}^{d}%
g(t,e_{i})z_{i}$ for almost every $t\in \lbrack0,T]$, where $z_{i}$ is the
$i$-th component of $z$.
\end{theorem}

For proving this theorem, we need the following lemmas. The first
lemma is a direct consequence of Jiang \cite{J} (see also
\cite{BCHMP, CCD, J1}).

\begin{lemma}\label{le5}
Suppose $g$ satisfies (H1)-(H3). If $\mathcal{E}_{g}[\xi
]=\mathcal{C}_{g}[\xi]$ for each $\xi \in L^{2}(\mathcal{F}_{T})$,
then $g$ is independent of $y$ and is positively homogeneous in $z.$
\end{lemma}

\begin{proof}
Since $\mathcal{E}_{g}=\mathcal{C}_{g},$ we have%
\[
\mathcal{E}_{g}[\xi+c]=\mathcal{E}_{g}[\xi]+c\quad \mbox{for
each}\ c\in \mathbb{R}%
;\  \mathcal{E}_{g}[\lambda \xi]=\lambda \mathcal{E}_{g}[\xi]\quad
\mbox{for each}\  \lambda \geq0.
\]
From this, we obtain the result (see Theorems 3.1 and 3.4 in Jiang
\cite{J}). The proof is complete.
\end{proof}

The next lemma is a property of BSDE (see \cite{P1}).

\begin{lemma}\label{le6}
Suppose $g$ satisfies (H1)-(H3). Let $\xi$ be a $k_{1}$-dimensional
$\mathcal{F}_{t_{0}}$-measurable random vector and $\eta$ be a $k_{2}%
$-dimensional $\mathcal{F}_{T}$-measurable random vector, where $t_{0}%
\in \lbrack0,T)$ and $k_{1},k_{2}\in \mathbb{N}$. Then for each $f\in
C_{b}(\mathbb{R}^{k_{1}}\times \mathbb{R}^{k_{2}})$, we have
\[
\mathcal{E}_{g}[f(\xi,\eta)|\mathcal{F}_{t}]=\mathcal{E}_{g}[f(x,\eta
)|\mathcal{F}_{t}]|_{x=\xi},\ t\in \lbrack t_{0},T].
\]

\end{lemma}

\begin{proof}
We outline the proof for the convenience of the reader. The proof is divided
into two steps.

Step 1: Let $\xi$ be simple random vector, i.e., $\xi=\sum_{i=1}^{N}%
x_{i}I_{A_{i}}$, where $\{x_{i}\}_{i=1}^{N}\subset
\mathbb{R}^{k_{1}}$ and $\{A_{i}\}_{i=1}^{N}$ is an
$\mathcal{F}_{t_{0}}$-partition of $\Omega$. Let
$(y_{t}^{i},z_{t}^{i})_{t\in \lbrack0,T]}$ denote the solution of
BSDE (\ref{e1}) corresponding to the generator $g$ and terminal
value $f(x_{i},\eta)$ with
$i=1,\ldots,N$. Then it is easy to verify that $(\sum_{i=1}^{N}y_{t}%
^{i}I_{A_{i}},\sum_{i=1}^{N}z_{t}^{i}I_{A_{i}})_{t\in \lbrack
t_{0},T]}$ is the solution of BSDE (\ref{e1}) on $[t_{0},T]$
corresponding to the generator $g$ and terminal value
$\sum_{i=1}^{N}f(x_{i},\eta)I_{A_{i}}$. Noting that
$f(\sum_{i=1}^{N}x_{i}I_{A_{i}},\eta)=\sum_{i=1}^{N}f(x_{i},\eta)I_{A_{i}}$,
then for $t\in \lbrack t_{0},T]$, we have%
\[
\mathcal{E}_{g}[f(\xi,\eta)|\mathcal{F}_{t}]=\sum_{i=1}^{N}\mathcal{E}%
_{g}[f(x_{i},\eta)|\mathcal{F}_{t}]I_{A_{i}}=\mathcal{E}_{g}[f(x,\eta
)|\mathcal{F}_{t}]|_{x=\xi}.
\]

Step 2: For general $\xi$, we can choose some simple random vectors $\xi
_{n}\rightarrow \xi$. Since $f\in C_{b}(\mathbb{R}^{k_{1}}\times \mathbb{R}%
^{k_{2}})$, by Lemma \ref{le1}, we get for $t\in \lbrack t_{0},T]$,%
\[
P-a.s.,\  \mathcal{E}_{g}[f(\xi_{n},\eta)|\mathcal{F}_{t}]\rightarrow
\mathcal{E}_{g}[f(\xi,\eta)|\mathcal{F}_{t}],\mathcal{E}_{g}[f(x,\eta
)|\mathcal{F}_{t}]|_{x=\xi_{n}}\rightarrow \mathcal{E}_{g}[f(x,\eta
)|\mathcal{F}_{t}]|_{x=\xi}.
\]
Thus $\mathcal{E}_{g}[f(\xi,\eta)|\mathcal{F}_{t}]=\mathcal{E}_{g}%
[f(x,\eta)|\mathcal{F}_{t}]|_{x=\xi}$. The proof is complete.
\end{proof}

\begin{remark}\label{re7}
Let $f_{n}\in C_{b}(\mathbb{R}^{k_{1}}\times \mathbb{R}^{k_{2}})$ be
uniformly bounded such that $f_{n}\rightarrow f$. Then by Lemma
\ref{le1}, we can easily prove that Lemma \ref{le6} still holds for
$f$.
\end{remark}

The following lemma plays an important role in proving the main theorem with
$d=1.$

\begin{lemma}\label{le8}
Suppose that $d=1$ and $g$ satisfies (H1)-(H3). If $\mathcal{E}_{g}%
[\xi]=\mathcal{C}_{g}[\xi]$ for each $\xi \in L^{2}(\mathcal{F}_{T})$, then for
each $t\in \lbrack0,T],$ $n\in \mathbb{N}$, we have%
\[
\mathcal{E}_{g}[I_{[W_{T}\geq-n]}+I_{[0\geq W_{T}\geq-n]}|\mathcal{F}%
_{t}]=\mathcal{E}_{g}[I_{[W_{T}\geq-n]}|\mathcal{F}_{t}]+\mathcal{E}%
_{g}[I_{[0\geq W_{T}\geq-n]}|\mathcal{F}_{t}].
\]

\end{lemma}

\begin{proof}
Let $W_{t,T}$ denote $W_{T}-W_{t}$. For each $a<b,$ it is easy to verify that
$I_{[W_{t,T}\geq a]}$ and $I_{[b\geq W_{t,T}\geq a]}$ are comonotonic. Then,
by $\mathcal{E}_{g}=\mathcal{C}_{g}$ and the comonotonic additivity of the
Choquet expectation, we have%
\begin{equation}\label{e2}
\mathcal{E}_{g}[I_{[W_{t,T}\geq a]}+I_{[b\geq W_{t,T}\geq a]}]=\mathcal{E}%
_{g}[I_{[W_{t,T}\geq a]}]+\mathcal{E}_{g}[I_{[b\geq W_{t,T}\geq a]}].
\end{equation}
On the other hand, for each $l_{1},l_{2}\in \mathbb{R},$ it is easy to show
that $f(x,y):=l_{1}I_{[x+y\geq-n]}+l_{2}I_{[0\geq x+y\geq-n]}$ satisfies the
condition in Remark \ref{re7}. Hence, we have%
\begin{equation}\label{e3}
\mathcal{E}_{g}[l_{1}I_{[W_{T}\geq-n]}+l_{2}I_{[0\geq W_{T}\geq-n]}%
|\mathcal{F}_{t}]=\mathcal{E}_{g}[l_{1}I_{[W_{t,T}\geq-n-\bar{a}]}%
+l_{2}I_{[-\bar{a}\geq W_{t,T}\geq-n-\bar{a}]}]|_{\bar{a}=W_{t}}.
\end{equation}
Combining (\ref{e3}) with (\ref{e2}) yields the result, and the
proof is complete.
\end{proof}

The following lemma is our main theorem with $d=1,$ which is an
extension of Lemma \ref{le3}.

\begin{lemma}\label{le9}
Suppose that $d=1$ and $g$ satisfies (H1)-(H3). Then $\mathcal{E}%
_{g}[\xi]=\mathcal{C}_{g}[\xi]$ for each $\xi \in L^{2}(\mathcal{F}_{T})$ if
and only if $g$ is independent of $y$ and is linear in $z$, i.e.,
$g(t,z)=g(t,1)z$ for almost every $t\in \lbrack0,T].$
\end{lemma}

\begin{proof}
If $g(t,z)=g(t,1)z,$ by the Girsanov Theorem, the $g$-expectation is a linear
mathematical expectation. Therefore, $\mathcal{E}_{g}[\xi]=\mathcal{C}_{g}%
[\xi]$ for each $\xi \in L^{2}(\mathcal{F}_{T})$ and the proof of
sufficient condition is complete. Now we prove the necessary
condition. By Lemma \ref{le5}, $g$ is independent of $y$. For each
$n\in \mathbb{N}$, consider the following
BSDEs:%
\begin{align*}
y_{t}^{n}  & =I_{[W_{T}\geq-n]}+I_{[0\geq W_{T}\geq-n]}+\int_{t}^{T}%
g(s,z_{s}^{n})ds-\int_{t}^{T}z_{s}^{n}dW_{s},\\
\tilde{y}_{t}^{n}  & =I_{[W_{T}\geq-n]}+\int_{t}^{T}g(s,\tilde{z}_{s}%
^{n})ds-\int_{t}^{T}\tilde{z}_{s}^{n}dW_{s},\\
\hat{y}_{t}^{n}  & =I_{[0\geq W_{T}\geq-n]}+\int_{t}^{T}g(s,\hat{z}_{s}%
^{n})ds-\int_{t}^{T}\hat{z}_{s}^{n}dW_{s}.
\end{align*}
By Lemma \ref{le8}, we have
$y_{t}^{n}=\tilde{y}_{t}^{n}+\hat{y}_{t}^{n}$ for each
$t\in \lbrack0,T].$ Form this, we have%
\begin{equation}\label{e4}
dP\times dt-a.s.,\quad g(t,\tilde{z}_{t}^{n}+\hat{z}_{t}^{n})=g(t,\tilde{z}%
_{t}^{n})+g(t,\hat{z}_{t}^{n}).
\end{equation}
On the other hand, it follows from Lemma \ref{le5} that $g$ is
positively homogeneous.
Thus we have for almost every $t\in \lbrack0,T],$%
\begin{equation}\label{e5}
g(t,z)=g(t,1)z^{+}+g(t,-1)z^{-},
\end{equation}
where $z^{+}=\max \{z,0\},$ $z^{-}=(-z)^{+}.$ Set $h(t):=g(t,1)+g(t,-1),$ by
(\ref{e4}) and (\ref{e5}), we have%
\begin{equation}\label{e6}
dP\times dt-a.s.,\quad
h(t)(\tilde{z}_{t}^{n}+\hat{z}_{t}^{n})^{-}=h(t)(\tilde
{z}_{t}^{n})^{-}+h(t)(\hat{z}_{t}^{n})^{-}.
\end{equation}
Also, $dP\times dt-a.s.$, $\tilde{z}_{t}^{n}=\frac{1}{\sqrt{2\pi(T-t)}}%
\exp(-\frac{(n+W_{t}+\int_{t}^{T}g(s,1)ds)^{2}}{2(T-t)})>0$ (see Lemma 8 in
\cite[Chen et al. (2005a)]{CCD}). This with (\ref{e6}) implies%
\begin{equation}\label{e7}
dP\times dt-a.s.,\quad h(t)I_{[\hat{z}_{t}^{n}<0]}=0.
\end{equation}
Let $(\bar{y}_{t},\bar{z}_{t})_{t\in \lbrack0,T]}$ denote the
solution of BSDE (\ref{e1}) corresponding to the generator $g$ and
terminal value $I_{[W_{T}\leq0]}.$ It follows from Lemma \ref{le1}
that $\hat{z}_{t}^{n}\rightarrow \bar{z}_{t}$ as $n\rightarrow
\infty$ in $L^{2}(0,T;\mathbb{R}).$ Thus we can choose
$n_{i}\rightarrow \infty$ such that $dP\times dt-a.s.$, $\hat{z}_{t}^{n_{i}%
}\rightarrow \bar{z}_{t}.$ Noting that $\bar{z}_{t}=\frac{-1}{\sqrt{2\pi(T-t)}%
}\exp(-\frac{(W_{t}-\int_{t}^{T}g(s,-1)ds)^{2}}{2(T-t)})<0,$ then by
(\ref{e7}), we can deduce that for almost every $t\in \lbrack0,T],$
$h(t)=0.$ Thus $g(t,z)=g(t,1)z$ for almost every $t\in \lbrack0,T].$
The proof is complete.
\end{proof}

\begin{corollary}\label{co10}
Suppose $g$ satisfies (H1)-(H3). If $\mathcal{E}_{g}[\xi]=\mathcal{C}_{g}%
[\xi]$ for each $\xi \in L^{2}(\mathcal{F}_{T})$, then $g$ is independent of
$y$ and is homogeneous in $z$, i.e., for almost every $t\in \lbrack0,T],$
$g(t,\lambda z)=\lambda g(t,z)$ for each $\lambda \in \mathbb{R}$.
\end{corollary}

\begin{proof}
For each fixed $a\in \mathbb{R}^{d}$ with $|a|=1$, set $\tilde{W}_{t}:=a\cdot
W_{t}$ and $\mathcal{\tilde{F}}_{t}:=\sigma \{ \tilde{W}_{s}:s\leq t\}$ for each
$t\in \lbrack0,T]$. Obviously, $(\tilde{W}_{t})_{t\in \lbrack0,T]}$ is a
$1$-dimensional Brownian motion. Define $\tilde{g}:[0,T]\times \mathbb{R\times
R\rightarrow R}$ by $\tilde{g}(t,y,z):=g(t,y,az)$. It is easy to verify that
$\tilde{g}$ satisfies (H1)-(H3). For each $\xi \in L^{2}(\mathcal{\tilde{F}%
}_{T})$, let $(y_{t},z_{t})_{t\in \lbrack0,T]}$ denote the solution of the
following BSDE:%
\[
y_{t}=\xi+\int_{t}^{T}\tilde{g}(s,y_{s},z_{s})ds-\int_{t}^{T}z_{s}d\tilde
{W}_{s},\ 0\leq t\leq T.
\]
Then it is easy to check that $(y_{t},az_{t})_{t\in \lbrack0,T]}$ is
the solution of BSDE (\ref{e1}) corresponding to the generator $g$
and terminal value
$\xi$. From this, we deduce that $\mathcal{E}_{g}[\xi]=\mathcal{E}_{\tilde{g}%
}[\xi]$ for each $\xi \in L^{2}(\mathcal{\tilde{F}}_{T})$. Noting that
$\mathcal{E}_{g}=\mathcal{C}_{g}$, we then get $\mathcal{E}_{\tilde{g}}%
[\xi]=\mathcal{C}_{\tilde{g}}[\xi]$ for each $\xi \in L^{2}(\mathcal{\tilde{F}%
}_{T})$. By Lemma \ref{le9}, we obtain
$\tilde{g}(t,y,z)=\tilde{g}(t,0,1)z$ for almost every $t\in
\lbrack0,T]$. Hence, by the Lipschitz assumption (H2), we have for
almost every $t\in \lbrack0,T]$, $g(t,y,\lambda a)=\lambda g(t,0,a)$
for each $\lambda \in \mathbb{R}$ and $a\in \mathbb{R}^{d}$ with
$|a|=1$, which implies that $g$ is independent of $y$ and is
homogeneous in $z$. We complete the proof.
\end{proof}

To prove the main theorem, we also pay more attention to the
following two lemmas.

\begin{lemma}\label{le11}
Suppose that $d=2$ and $g$ satisfies (H1)-(H3). If $\mathcal{E}_{g}%
[\xi]=\mathcal{C}_{g}[\xi]$ for each $\xi \in L^{2}(\mathcal{F}_{T})$, then for
each $\lambda \in \lbrack0,1]$, $t\in \lbrack0,T],$ $n\in \mathbb{N},$ we have%
\[
\mathcal{E}_{g}[I_{[W_{T}^{1}\geq n]}+\lambda I_{[W_{T}^{2}\geq0]}%
|\mathcal{F}_{t}]=\lambda \mathcal{E}_{g}[I_{[W_{T}^{1}\geq n]}+I_{[W_{T}%
^{2}\geq0]}|\mathcal{F}_{t}]+(1-\lambda)\mathcal{E}_{g}[I_{[W_{T}^{1}\geq
n]}|\mathcal{F}_{t}],
\]
where $W_{t}^{i}$ is the $i$-th component of $W_{t}$ with $i=1,2.$
\end{lemma}

\begin{proof}
Let $W_{t,T}^{i}$ denote $W_{T}^{i}-W_{t}^{i}$ with $i=1,2$. For
each fixed $\lambda \in \lbrack0,1],$ $a,b\in \mathbb{R}$, it is
easy to check that $(1-\lambda)I_{[W_{t,T}^{1}\geq a]}$ and
$\lambda(I_{[W_{t,T}^{1}\geq a]}+I_{[W_{t,T}^{2}\geq b]})$ are
comonotonic. The rest of the proof runs as in Lemma \ref{le8}, and
the proof is complete.
\end{proof}

\begin{lemma}\label{le12}
Suppose that $d=2$ and $g$ satisfies (H1)-(H3). If $\mathcal{E}_{g}%
[\xi]=\mathcal{C}_{g}[\xi]$ for each $\xi \in L^{2}(\mathcal{F}_{T})$, then $g$
is independent of $y$ and is linear in $z$, i.e., $g(t,z_{1},z_{2}%
)=g(t,1,0)z_{1}+g(t,0,1)z_{2}$ for almost every $t\in \lbrack0,T].$
\end{lemma}

\begin{proof}
It follows from Lemma \ref{le5} that $g$ is independent of $y$. For
each fixed
$\lambda \in(0,1),$ $n\in \mathbb{N},$ consider the following BSDEs:%
\begin{align*}
y_{t}^{\lambda,n} &  =I_{[W_{T}^{1}\geq n]}+\lambda I_{[W_{T}^{2}\geq0]}%
+\int_{t}^{T}g(s,z_{1,s}^{\lambda,n},z_{2,s}^{\lambda,n})ds-\int_{t}%
^{T}z_{1,s}^{\lambda,n}dW_{s}^{1}-\int_{t}^{T}z_{2,s}^{\lambda,n}dW_{s}^{2},\\
\tilde{y}_{t}^{n} &  =I_{[W_{T}^{1}\geq n]}+I_{[W_{T}^{2}\geq0]}+\int_{t}%
^{T}g(s,\tilde{z}_{1,s}^{n},\tilde{z}_{2,s}^{n})ds-\int_{t}^{T}\tilde{z}%
_{1,s}^{n}dW_{s}^{1}-\int_{t}^{T}\tilde{z}_{2,s}^{n}dW_{s}^{2},\\
\hat{y}_{t}^{n} &  =I_{[W_{T}^{1}\geq n]}+\int_{t}^{T}g(s,\hat{z}_{1,s}%
^{n},\hat{z}_{2,s}^{n})ds-\int_{t}^{T}\hat{z}_{1,s}^{n}dW_{s}^{1}-\int_{t}%
^{T}\hat{z}_{2,s}^{n}dW_{s}^{2}.
\end{align*}
By Lemma \ref{le11}, we have $y_{t}^{\lambda,n}=\lambda \tilde{y}_{t}^{n}%
+(1-\lambda)\hat{y}_{t}^{n}$ for each $t\in \lbrack0,T]$. From this, we deduce
that $dP\times dt-a.s.$,
\[
g(t,\lambda \tilde{z}_{1,t}^{n}+(1-\lambda)\hat{z}_{1,t}^{n},\lambda \tilde
{z}_{2,t}^{n}+(1-\lambda)\hat{z}_{2,t}^{n})=\lambda g(t,\tilde{z}_{1,t}%
^{n},\tilde{z}_{2,t}^{n})+(1-\lambda)g(t,\hat{z}_{1,t}^{n},\hat{z}_{2,t}^{n}).
\]
Since $\lambda \in(0,1)$ is arbitrary, by Lemma \ref{le5}, we obtain
that
$dP\times dt-a.s.$,%
\begin{equation}\label{e8}
g(t,\tilde{z}_{1,t}^{n}+l\hat{z}_{1,t}^{n},\tilde{z}_{2,t}^{n}+l\hat{z}%
_{2,t}^{n})=g(t,\tilde{z}_{1,t}^{n},\tilde{z}_{2,t}^{n})+g(t,l\hat{z}%
_{1,t}^{n},l\hat{z}_{2,t}^{n})\ \mbox{for each}\ l\geq0.
\end{equation}
Noting that $g(t,z_{1},0)=g(t,1,0)z_{1}$ for almost every $t\in \lbrack0,T],$
then we have%
\begin{equation}\label{e9}
dP\times dt-a.s.,\ (\hat{z}_{1,t}^{n},\hat{z}_{2,t}^{n})=(\frac{1}{\sqrt
{2\pi(T-t)}}\exp(-\frac{(n-W_{t}^{1}-\int_{t}^{T}g(s,1,0)ds)^{2}}{2(T-t)}),0).
\end{equation}
Combining (\ref{e8}) with (\ref{e9}), we get
\begin{equation}\label{e10}
dP\times dt-a.s.,\ g(t,\tilde{z}_{1,t}^{n}+p,\tilde{z}_{2,t}^{n}%
)=g(t,\tilde{z}_{1,t}^{n},\tilde{z}_{2,t}^{n})+g(t,p,0)\ \mbox{for
each}\ p\geq0.
\end{equation}
Let $(\bar{y}_{t},\bar{z}_{1,t},\bar{z}_{2,t})_{t\in \lbrack0,T]}$
be the solution of BSDE (\ref{e1}) corresponding to the generator
$g$ and terminal value $I_{[W_{T}^{2}\geq0]}$. By Lemma \ref{le1},
we have $(\tilde{z}_{1,t}^{n},\tilde
{z}_{2,t}^{n})\rightarrow(\bar{z}_{1,t},\bar{z}_{2,t})$ in $L^{2}%
(0,T;\mathbb{R}^{2}).$ Since $g$ satisfies Lipschitz assumption (H2), we get
for each $p\geq0,$%
\[
g(t,\tilde{z}_{1,t}^{n}+p,\tilde{z}_{2,t}^{n})\rightarrow g(t,\bar{z}%
_{1,t}+p,\bar{z}_{2,t})\ \mbox{in}\ L^{2}(0,T;\mathbb{R}).
\]
This with (\ref{e10}) implies that
\begin{equation}\label{e11}
dP\times dt-a.s.,\ g(t,\bar{z}_{1,t}+p,\bar{z}_{2,t})=g(t,\bar{z}_{1,t}%
,\bar{z}_{2,t})+g(t,p,0)\ \mbox{for each}\ p\geq0.
\end{equation}
Also, we have%
\begin{equation}\label{e12}
dP\times dt-a.s.,\ (\bar{z}_{1,t},\bar{z}_{2,t})=(0,\frac{1}{\sqrt{2\pi(T-t)}%
}\exp(-\frac{(W_{t}^{2}+\int_{t}^{T}g(s,0,1)ds)^{2}}{2(T-t)})).
\end{equation}
It follows from (\ref{e11}), (\ref{e12}) and Lemma \ref{le5} that
for almost every $t\in
\lbrack0,T],$%
\begin{equation}\label{e13}
g(t,p,1)=g(t,0,1)+g(t,p,0)\ \mbox{for each}\ p\geq0.
\end{equation}
From (\ref{e13}) and Corollary \ref{co10}, we can easily deduce that
for almost every
$t\in \lbrack0,T],$%
\begin{equation}\label{e14}
g(t,z_{1},z_{2})=g(t,1,0)z_{1}+g(t,0,1)z_{2}\ \mbox{for
each}\ z_{1}\cdot z_{2}%
\geq0.
\end{equation}
On the other hand, set $\tilde{W}_{t}:=(W_{t}^{1},-W_{t}^{2})$ and $\tilde
{g}(t,z_{1},z_{2})=g(t,z_{1},-z_{2}).$ Analysis similar to that in the proof
of Corollary \ref{co10} shows that $\mathcal{E}_{\tilde{g}}[\xi]=\mathcal{C}%
_{\tilde{g}}[\xi]$ for each $\xi \in L^{2}(\mathcal{F}_{T})$. Then
we have (\ref{e14})
for $\tilde{g}$, which gives that for almost every $t\in \lbrack0,T],$%
\[
g(t,z_{1},z_{2})=g(t,1,0)z_{1}+g(t,0,1)z_{2}\ \mbox{for
each}\ z_{1}\cdot z_{2}%
\leq0.
\]
The proof is now complete.
\end{proof}

We now prove the main theorem.

\textbf{Proof of Theorem \ref{th4}. }The sufficient condition can be
proved by the same method as in Lemma \ref{le9}. We only prove the
necessary condition. For $d=2,$ by Lemma \ref{le12}, the result
holds. We only prove the case $d>2.$ For each fixed
$a\in \mathbb{R}^{d-1}$ with $|a|=1,$ set $\tilde{W}_{t}:=(a\cdot(W_{t}%
^{1},\ldots,W_{t}^{d-1}),W_{t}^{d})$ and $\mathcal{\tilde{F}}_{t}%
:=\sigma \{ \tilde{W}_{s}:s\leq t\}$ for each $t\in \lbrack0,T]$. By
Lemma \ref{le5}, $g$ is independent of $y$, we define
\thinspace$\tilde{g}:[0,T]\times
\mathbb{R\times R\rightarrow R}$ by $\tilde{g}(t,z_{1},z_{2}):=g(t,az_{1}%
,z_{2})$. As in the proof of Corollary \ref{co10}, we can get $\mathcal{E}_{\tilde{g}%
}[\xi]=\mathcal{C}_{\tilde{g}}[\xi]$ for each $\xi \in
L^{2}(\mathcal{\tilde
{F}}_{T})$. By Lemma \ref{le12}, we have for almost every $t\in \lbrack0,T],$%
\[
\tilde{g}(t,z_{1},z_{2})=\tilde{g}(t,1,0)z_{1}+\tilde{g}(t,0,1)z_{2}.
\]
Since $a$ is arbitrary, by Corollary \ref{co10}, we obtain for
almost every
$t\in \lbrack0,T],$%
\[
g(t,z_{1},\ldots,z_{d-1},z_{d})=g(t,z_{1},\ldots,z_{d-1},0)+g(t,e_{d})z_{d}.
\]
Define $\bar{g}:[0,T]\times \mathbb{R}^{d-1}\mathbb{\rightarrow R}$ by $\bar
{g}(t,z):=g(t,z,0).$ We now apply the above argument again, with $g$ replaced
by $\bar{g}$, to obtain that for almost every $t\in \lbrack0,T],$%
\[
\bar{g}(t,z_{1},\ldots,z_{d-2},z_{d-1})=\bar{g}(t,z_{1},\ldots,z_{d-2}%
,0)+\bar{g}(t,0,\ldots,0,1)z_{d-1},
\]
that is%
\[
g(t,z_{1},\ldots,z_{d-2},z_{d-1},0)=g(t,z_{1},\ldots,z_{d-2},0,0)+g(t,e_{d-1}%
)z_{d-1}.
\]
Continuing this process, we can prove that $g(t,z)=\sum_{i=1}^{d}%
g(t,e_{i})z_{i}$ for almost every $t\in \lbrack0,T].$ The proof is complete.
$\Box$

\end{document}